\documentclass[12pt]{article}
\usepackage{amsfonts, amsmath, amstext, amssymb}

\usepackage[dvips]{graphicx}
\usepackage{latexsym}

\newtheorem{thm}{Theorem}[section]
\newtheorem{prop}[thm]{Proposition}

\newtheorem{cor}[thm]{Corollary}

\newtheorem{conj}[thm]{Conjecture}

\newtheorem{definitiontemp}[thm]{Definition}
\newenvironment{defn}{\begin{definitiontemp}
\normalfont}{\end{definitiontemp}}

\newcommand{\urd}{uniformly relatively decidable} 

\newcommand{\proves}{\vdash}

\newcommand{\Th}[1]{\text{Th}(#1)}

\newenvironment{pf}{\begin{trivlist}\item[\hskip\labelsep
{\it Proof.}]}{\end{trivlist}}

\newcommand{\set}[2]{\ensuremath{ \{ #1 : #2 \} }}

\newcommand{\Z}{\mathbb{Z}}

\newcommand{\Q}{\mathbb{Q}}
\newcommand{\R}{\mathcal{R}}
\newcommand{\A}{\mathcal{A}}
\newcommand{\B}{\mathcal{B}}
\newcommand{\C}{\mathcal{C}}
\newcommand{\D}{\mathcal{D}}

\renewcommand{\P}{\mathcal{P}}

\renewcommand{\L}{\mathcal{L}}

\newcommand{\avec}{\vec{a}}
\newcommand{\bvec}{\vec{b}}
\newcommand{\cvec}{\vec{c}}
\newcommand{\dvec}{\vec{d}}

\newcommand{\gvec}{\vec{g}}

\newcommand{\xvec}{\vec{x}}

\newcommand{\yvec}{\vec{y}}

\newcommand{\zvec}{\vec{z}}

\newcommand{\la}{\langle}
\newcommand{\ra}{\rangle}

\newcommand{\Adj}{\textbf{Adj}}

\def\diverges{\!\uparrow}
\def\converges{\!\downarrow}
\newcommand{\at}{\char'100}

\def\frakS{{\mathfrak S}}

\newcommand{\qed}{\hbox to 0pt{}\nobreak\hfill\rule{2mm}{2mm}}

\newcommand{\comment}[1]{}

\newcommand{\dom}[1]{\text{dom}(#1)}
\newcommand{\rg}[1]{\text{range}(#1)}

\newcommand{\bfc}{\boldsymbol{c}}
\newcommand{\bfd}{\boldsymbol{d}}

\def\s01{\ensuremath{\Sigma^0_1}}
\def\d02{\ensuremath{\Delta^0_2}}
\def\phi{\varphi}

\def\res{\!\!\upharpoonright\!}

\title{Model Completeness and Relative Decidability}
\author{Jennifer Chubb, Russell Miller\thanks{The second author was supported
by NSF grant \# DMS-1362206, Simons Foundation grant \# 581896,
and several PSC-CUNY research awards.}, \& Reed Solomon}

\begin{document}

\maketitle

\begin{abstract}
We study the implications of model completeness of a theory
for the effectiveness of presentations of models of that theory.
It is immediate that for a computable model $\A$ of a computably enumerable,
model complete theory, the entire elementary diagram $E(\A)$ must be decidable.
We prove that indeed a c.e.\ theory $T$ is model complete if and only if
there is a uniform procedure that succeeds in deciding $E(\A)$ from the
atomic diagram $\Delta(\A)$ for all countable models $\A$ of $T$.  Moreover, if
every presentation of a single isomorphism type $\A$ has this property
of relative decidability, then there must be a procedure with succeeds
uniformly for all presentations of an expansion $(\A,\avec)$ by finitely
many new constants.  We end with a conjecture about the situation
when all models of a theory are relatively decidable.
\end{abstract}

\section{Introduction}
\label{sec:intro}

The broad goal of computable model theory is to investigate the
effective aspects of model theory.  Here we will carry out exactly
this process with the model-theoretic concept of \emph{model completeness}.
This notion is well-known and has been widely studied in model theory,
but to our knowledge there has never been any thorough examination
of its implications for computability in structures with the domain $\omega$.
We now rectify this omission, and find natural and satisfactory equivalents
for the basic notion of model completeness of a first-order theory.
Our two principal results are each readily stated:  that a computably enumerable
theory is model complete if and only if there is a uniform procedure for
deciding the elementary diagram of each of its models (on the domain $\omega$)
from the corresponding atomic diagram (Theorem \ref{thm:infQE});
and that if every presentation of a particular isomorphism type $\A$
individually has this property (which we call \emph{relative decidability}),
then for some expansion $(\A,\avec)$ by finitely many constants, there
must be a uniform version of the procedure that succeeds for all presentations
(Theorem \ref{thm:URD}).

Recall that a structure $\A$ with domain $\omega$ (or a decidable subset of $\omega$)
is \emph{computable} if its atomic diagram is computable, and is \emph{decidable}
if its elementary diagram is computable.  This choice of terminology may seem arbitrary,
but it is well established, and we will maintain it here.  In our notation,
$\Delta(A)$ is the atomic diagram of a structure $\A$ with domain $\omega$.
In a fixed countable language $\L$, fix an effective G\"odel coding of all atomic sentences
in the language $(\L,c_0,c_1,\ldots)$ with new constants for the naturals; then $\Delta(\A)$
can be viewed as a subset of $\omega$, containing the code numbers of those
sentences true in $\A$, with each $c_n$ representing the domain element $n$ of $\A$.
We assume implicitly that $\Delta(\A)$ always respects equality on $\omega$,
i.e., that the atomic statement $c_i=c_j$ lies in $\Delta(\A)$ only if $i=j$.
$E(\A)$ is the elementary diagram (a.k.a.\ the complete diagram) of $\A$,
coded similarly as a subset of $\omega$.

\begin{defn}
\label{defn:QE}
A theory $T$ \emph{has quantifier elimination down to $\C$}, where $\C$
is a class of formulas, if, for every (finitary) formula $\alpha(\xvec)$,
there exists a formula $\gamma(\xvec)$ in $\C$, with the same free variables
$\xvec$, such that
$$ T\models \forall\xvec (\alpha(\xvec)\leftrightarrow \gamma(\xvec)).$$
\end{defn}
For example, it is well-known that a theory is model complete
if and only if it has quantifier elimination down to (finitary)
existential formulas, or equivalently, down to universal formulas.
(For a proof, see \cite[Theorem 3.5.1]{CK90}.)  Later on, we will
consider quantifier elimination down to existential $L_{\omega_1\omega}$
formulas.

Having introduced the term, we recall the definition.
\begin{defn}
\label{defn:modelcomplete}
A theory $T$ is \emph{model complete} if, for every model $\B$ of $T$,
every substructure $\A\subseteq\B$ which is itself a model of $T$ is an
elementary substructure of $\B$.
\end{defn}
It is equivalent to require that, for every $\A\subseteq\B$ as described,
every existential statement $\exists\xvec\phi(\avec,\xvec)$ which holds
in $\B$ (of a tuple $\avec$ from $\A$) also holds in $\A$.

Our two principal theorems will be proven in Sections \ref{sec:basics}
and \ref{sec:thms}, respectively.  In between, in Section \ref{sec:examples},
we offer several examples of theories that are or are not model complete;
these might be perused before Section \ref{sec:basics} by the reader
desiring a refresher, but they also lead in well to Section \ref{sec:thms}.
At the end, in Section \ref{sec:precomplete}, we offer a conjecture
regarding the situation of a theory for which every model is relatively
decidable, but not uniformly so.  This conjecture involves a notion
called \emph{model precompleteness} that we believe to be new
in the literature, and we hope to encourage model theorists to examine it.

Arguments in this article always use countable structures, whether we say so or not.
A \emph{presentation} of a countable structure is simply a structure isomorphic
to it, whose domain is the set $\omega$ of nonnegative integers.
This enables our G\"odel coding above to make sense in all cases.
Finally, all languages considered here will be computable languages,
with at most countably many symbols.  We believe that, for a noncomputable
countable language $\L$, all the results would go through without difficulty
if one simply relativized the statements and arguments to the Turing degree of $\L$.

\section{Basics of Model Completeness}
\label{sec:basics}

We begin with the following result, which was remarked in \cite[Prop.\ 6.7]{H98}.

\begin{prop}
\label{prop:RD}
Let $T$ be a model complete theory, and assume $T$ is c.e.
Then every computable model of $T$ is decidable.
\end{prop}
\begin{pf}
Fix any formula $\phi(x_0,\ldots,x_{n-1})$ and any $\avec\in\omega^n$
for which we wish to determine whether $\phi(\avec)$ holds in $\A$.
As noted above, $T$ has the property that,
for this (and every) formula $\phi(y_0,\ldots,y_{n-1})$, there is a universal formula
in the same variables which is provably equivalent (under $T$) to $\phi$.
Since $T$ is c.e., we may therefore search until we find quantifier-free
formulas $\alpha(x_0,\ldots,x_{n-1},y_0,\ldots,y_m)$ and $\beta(\xvec,\yvec)$
and a deduction showing
$$ T \proves \forall\xvec[(\phi(\xvec)\leftrightarrow\forall\yvec~\alpha(\xvec,\yvec))~\&~
((\neg\phi(\xvec))\leftrightarrow\forall\yvec~\beta(\xvec,\yvec))].$$
Model-completeness of $T$ ensures that this search will terminate.

Now, in the computable structure $\A$, we search for some tuple
$\bvec\in\omega^m$ such that either $(\neg\alpha(\avec,\bvec))$
or $(\neg\beta(\avec,\bvec))$ holds in $\A$.  By our choice of $\alpha$ and $\beta$,
this search must also terminate.  When it does, we know which of
the formulas $\forall\yvec~\alpha(\avec,\yvec))$ and $\forall\yvec~\beta(\avec,\yvec))$
fails to hold in $\A$, and from this we determine which of $\phi(\avec)$
and $(\neg\phi(\avec))$ holds in $\A$.
\qed\end{pf}

The proposition generalizes easily to the following full theorem.
(It would be natural here for $\bfc$ to be the Turing degree
of an axiomatization of the theory $T$.)
\begin{thm}
\label{thm:forwards}
Let $T$ be a model complete theory, and $\bfc$ a Turing degree.
Assume $T$ is $\bfc$-c.e.
Then, for every Turing degree $\bfd\geq\bfc$, every $\bfd$-computable
model of $T$ is $\bfd$-decidable.
\qed\end{thm}

We ask whether the converse of this theorem holds:  if a c.e.\ theory $T$
has the property that every $\bfd$-computable model of $T$ is $\bfd$-decidable,
must $T$ be model complete?  (Likewise for $\bfc$-c.e.\ theories where this holds for all $\bfd\geq\bfc$.)
Alternatively, there might be some version involving $L_{\omega_1\omega}$
formulas, perhaps computable ones, which would be equivalent to model completeness.

The first answer is that the converse does not hold.
\begin{prop}
\label{prop:noconverse}
The theory $T=\Th{\omega, S}$ is not model complete, yet for every $\bfd$,
every $\bfd$-computable model is $\bfd$-decidable.
\end{prop}
\begin{pf}
Note that the theory $T$ is computable because $\Th{\omega, S, 0}$ is a computable theory which 
is a conservative extension of $T$. Therefore, for any $\{ S \}$-formula 
$\varphi$, $\varphi \in T$ if and only if $\varphi \in \Th{\omega, S, 0}$. 
 
If $\A=(\omega,S)$ and $\B=(\omega-\{0\},S)$, then $\A$ and $\B$ are both
models of $T$, with $\A\subseteq\B$, yet $\B$ is not an elementary substructure
of $\A$:  the element $1$ is a successor in $\A$, but not in $\B$. Therefore, $T$ is not model complete.

However, whenever a model $\C$ of $T$ is $\bfd$-computable, there is a decision
procedure for its elementary diagram.  Indeed, all models of the model complete
theory $\Th{\omega,0,S}$ have this property, and one simply applies the same procedure
to $\C$ using the unique element of $\C$ with no predecessor as the constant.
\qed\end{pf}

This proposition points to the true kernel of the converse.  The process of deciding the
elementary diagram of a model of $\Th{\omega,S}$ from its atomic diagram
was nonuniform:  it required knowledge of the unique non-successor element,
which cannot be determined uniformly from the atomic diagram.  On the other hand,
the procedure described in Theorem \ref{thm:forwards} was uniform.
Therefore, we modify our conjecture about the converse to require uniformity.

\begin{defn}
\label{defn:urd}
A structure $\A$ is \emph{relatively decidable} if $E(\A)\leq_T\Delta(\A)$.
A class $\mathfrak{S}$ of structures is \emph{\urd} if there exists a single Turing
functional $\Gamma$ such that, for every $\A\in\frakS$, the function
$\Gamma^{\Delta(\A)}$ is the characteristic function of the elementary diagram $E(\A)$.

The structure $\A$ is \emph{\urd} if the class of all structures (with domain $\omega$)
isomorphic to $\A$ is \urd.
A theory $T$ is \emph{\urd} if the class of all models of $T$
with domain $\omega$ is \urd.
\end{defn}

\comment{
\begin{conj}
A c.e.\ theory $T$ is model complete if and only if $T$ is \urd.
\end{conj}

The forwards direction is immediate; the converse is not.  To address the converse,
}
To prove equivalence between this concept and model completeness,
we use a broader form of quantifier elimination.

\begin{thm}
\label{thm:infQE}
For a computably enumerable theory $T$ (in particular, for any theory
with a decidable axiom set), the following are equivalent:
\begin{enumerate}
\item
$T$ is model complete.
\item
$T$ is \urd.
\item
$T$ has effective quantifier elimination down to
$\Sigma^c_1$ formulas, that is, down to computable infinitary $\Sigma_1$ formulas.
\item
$T$ has effective quantifier elimination down to finitary $\Sigma_1$ formulas.
\end{enumerate}
\end{thm}

\comment{
Notice that here we must assue the quantifier elimination to be an effective procedure,
as it no longer is reasonable to search for a proof from $T$ of the equivalence of
a given $\alpha$ to a $\Sigma^c_1$ formula $\beta$.  On the other hand, since
this search is now unnecessary, there is no longer any need for $T$ to be c.e.
}

\begin{pf}
(1) and (4) are well-known to be equivalent; see for instance \cite[Theorem 3.5.1]{CK90}.
(The effectiveness of the quantifier elimination is not stated there, but
a c.e.\ theory with quantifier elimination always has effective quantifier elimination.)
Moreover, the proof of Proposition \ref{prop:RD} actually shows that
(4) implies (2), since the procedure given there for deciding the elementary
diagram of $\A$ from $\Delta(\A)$ is uniform for all models $\A$ of $T$.

It is also quickly seen that (3) implies (1).  Recall Definition \ref{defn:modelcomplete},
and let $\A$ and $\B$ be models of $T$ with $\A\subseteq\B$.  Now for each
formula $\phi(\xvec)$, we have a $\Sigma_1^c$ formula $\gamma(\xvec)$
equivalent to $(\neg\phi(\xvec))$ in all models of $T$.  For each tuple $\avec$
from $\A$ such that $\B\models\phi(\avec)$, it is immediate that $\B\not\models\gamma(\avec)$,
hence that $\A\not\models\gamma(\avec)$ (because $\Sigma^1_c$ formulas true
in $\A$ must hold in superstructures of $\A$), hence that $\A\models\phi(\avec)$.
Thus $\A$ is an elementary substructure of $\B$.

\comment{
First, suppose that $T$ has the effective quantifier elimination described.
Then, given a model $\A$ of $T$ and an arbitrary finitary statement $\alpha(\avec)$
with $\avec\in\omega^n$, we use this procedure to determine $\Sigma^c_1$
formulas $\beta(\xvec)$ and $\gamma(\xvec)$ such that
$$ T\models \forall\xvec[(\alpha(\xvec)\leftrightarrow\beta(\xvec))~\&~
((\neg\alpha(\xvec))\leftrightarrow\gamma(\xvec))].$$
Using the $\Delta(\A)$-oracle, we can search for witnesses in $\A$
to either $\beta(\avec)$ or $\gamma(\avec)$; eventually this search
must terminate, telling us whether $\alpha(\avec)$ holds in $\A$ or not.
}

It remains to show that (2) implies (3).
Assume that $T$ is \urd, via some functional $\Gamma$.  
Consider the arbitrary (finitary) formula $\alpha(x_0,\ldots,x_n)$.
For simplicity, we will assume for the time being that $\alpha$
includes the conditions $x_i\neq x_j$ for all $i<j\leq n$.  At the end of the proof,
we show how to remove this assumption.

We claim
that in every model of $T$, the formula $\alpha(\xvec)$
will be equivalent to the following $\Sigma^c_1$ formula $\beta_{\alpha}(\xvec)$:
$$
\bigvee\!\!\!\!\!\!\!\!\!\!\bigvee_{\sigma\in H_{\alpha}}
\exists y_{n+1},\ldots,y_{m_\sigma}~\gamma_{\sigma}(\xvec,\yvec),$$
where the set $H_{\alpha}$ is c.e., uniformly in $\alpha$, as we now explain.
First, given any $\sigma\in 2^{<\omega}$, let $\gamma_\sigma(c_0,\ldots,c_{m_\sigma})$
be the finitary quantifier-free formula
$$ \left(\bigwedge_{\sigma(\ulcorner \psi\urcorner)=1}\psi\right)~\bigwedge~
\left(\bigwedge_{\sigma(\ulcorner\psi\urcorner)=0}(\neg\psi)\right)~\bigwedge~
\left(\bigwedge_{i<j<|\sigma|}c_i\neq c_j\right),$$
using our G\"odel numbering of the atomic formulas in the language of $T$
with new constants $c_0,c_1,\ldots$.  Here $\sigma$ is viewed as a possible
initial segment of an atomic diagram, although of course for many strings $\sigma$,
$\gamma_\sigma$ may already contradict $T$.
Next, let
$$ H_{\alpha} = \set{\sigma\in 2^{<\omega}}
{\Gamma^{\sigma}(\ulcorner \alpha(c_0,\ldots,c_n)\urcorner)\converges=1}.$$
So $\sigma$ lies in $H_\alpha$ if $\gamma_\sigma$ is enough information
for the procedure $\Gamma$ to conclude that $\alpha(\cvec)$ holds in those models
$\A$ of $T$ (if any) whose atomic diagrams specify that $\gamma_\sigma$ holds in $\A$.

Finally, when defining $\beta_{\alpha}(\xvec)$, we quantified over certain $y_i$;
we now explain how.  Each $\gamma_\sigma$ is a quantifier-free formula
involving some constants $c_0,\ldots,c_{m_\sigma}$, quite possibly with $m_\sigma>n$.
In $\beta_{\alpha}(\xvec)$, each remaining $c_i$ (with $i>n$)
is replaced by $y_i$, with an existential quantification $\exists y_i$ added
in front of $\gamma_\sigma$.  Thus the formula $\beta_{\alpha}(\xvec)$
involves $\xvec$ but no other free variables, nor any constants.

It is clear that $H_{\alpha}$ is a c.e.\ set, and therefore that $\beta_{\alpha}(\xvec)$
is a $\Sigma^c_1$ formula.  Moreover, the process is uniform in $\alpha$.
It remains to show that, for every model $\A$ of $T$
and every $\avec$ from $\A$, $\alpha(\avec)$
holds in $\A$ if and only $\beta_{\alpha}(\avec)$ holds there.
Suppose first that $\A\models\alpha(\avec)$.  
Fix a permutation $f$ of $\omega$ satisfying:
\begin{itemize}
\item
$f(a_i)=i$ for all $i\leq n$; and
\item
$f(a)=a$ for all but finitely many $a\in\omega$.
\end{itemize}
(Recall that we are assuming for the time being that $\alpha(\xvec)$
includes the conditions $x_i\neq x_j$ for $i\neq j$.  This is necessary for
$f$ to exist.)
Define the structure $\A_f$ so that $f$ is an isomorphism from $\A$
onto $\A_f$.  Then $\A_f \models \alpha(0,\ldots,n)$,
and therefore $\Gamma^{\Delta(\A_f)}(\ulcorner\alpha(0,\ldots,n)\urcorner)\converges=1$.
Let $\sigma$ be the initial segment of $\Delta(\A_f)$ as long as the use
of this computation, so that $\sigma\in H_\alpha$.  But this means that
$\A_f \models \gamma_\sigma(0,\ldots, m_\sigma)$, and so the
isomorphism $f^{-1}$ shows that $\A$ also satisfies
$\A \models \exists y_{n+1}\cdots\exists y_{m_\sigma} \gamma_\sigma(\avec,\yvec)$.
Thus $\A$ satisfies $\beta_\alpha(\avec)$, as required.

Conversely, suppose $\A \models \beta_{\alpha}(\avec)$.
Then for some $\sigma\in H_{\alpha}$, 
$\A \models \exists\yvec~\gamma_\sigma(\avec,\yvec)$. Fix $b_{n+1}, \ldots, b_{m_{\sigma}}$ such that 
$\A \models \gamma_{\sigma}(\avec,\bvec)$. Define $f$ to be a permutation as above with $f(a_i) = i$ and $f(b_j) = j$. Let 
$\A_f$ be the corresponding copy of $\A$, so $\A_f \models \gamma_{\sigma}(0, \ldots, m_{\sigma})$ and hence 
$\sigma$ is an initial segment of $\Delta(\A_f)$. Because $\sigma\in H_{\alpha}$,
we see that $\Gamma^{\Delta(\A_f)}(\ulcorner\alpha(0,\ldots,n)\urcorner)\converges=1$.
This means that $\alpha(0,\ldots,n)$ must hold in $\A_f$, and the isomorphism
$f^{-1}$ now shows that $\A \models \alpha(a_0,\ldots,a_n)$ as well.

All of this was proven under the assumption that $\alpha(\xvec)$ includes
the conditions $x_i\neq x_j$ for $i<j\leq n$.  For the general case, one simply
expresses $\alpha$ as a finite disjunction of the possibilities, and applies
the process above to each one individually, getting a finite disjunction of
$\Sigma^c_1$ formulas.  For example, if $n=2$, then
$\alpha(x_0,x_1,x_2)$ is equivalent to the following disjunction.
\begin{align*}
&\alpha(x_0,x_0,x_0)\\
\text{or }&(\alpha(x_0,x_1,x_1)~\&~x_0\neq x_1)\\
\text{or }&(\alpha(x_0,x_1,x_0)~\&~x_0\neq x_1)\\
\text{or }&(\alpha(x_0,x_0,x_1)~\&~x_0\neq x_1)\\
\text{or }
&(\alpha(x_0,x_1,x_2)~\&~x_0\neq x_1~\&~x_0\neq x_2~\&~x_1\neq x_2).~~~~~~~~~~~~~~
\qed
\end{align*}
\end{pf}

Theorem \ref{thm:infQE} has a natural generalization to theories $T$
that are not computably enumerable.  If $S$ is a set that can enumerate $T$
(for example, if $S$ is an axiom set for $T$), then $T$ is model complete
if and only if it is uniformly relatively $S$-decidable.  By this we mean
that there is a Turing functional $\Gamma$ such that, for every model
$\A$ of $T$ with domain $\omega$, $\Gamma^{S\oplus\Delta(\A)}$
computes $E(\A)$.  The argument is identical to that in the proof
of the theorem, and uses $S$-effective quantifier elimination down to
$S$-computable infinitary $\Sigma_1$ formulas as an intermediate
equivalent.

\section{Examples}
\label{sec:examples}

Here we provide some examples of situations where the addition
of new symbols to the language can allow a theory to become model complete.
Several examples of model complete theories appear in \cite[\S 3.5]{CK90}.
One of these is the dense linear order with endpoints, in a language with the relation
$<$ and two constant symbols to name the endpoints.  Without the constants,
the theory $T$ of dense linear orders with endpoints would not be model complete:
the rational interval $[0,\frac12]_{\Q}$ would be a substructure of $[0,1]_{\Q}$, and both
would be models of $T$, yet the larger structure would satisfy $(\exists x)~\frac12<x$
and the substructure would not.  When $T$ is augmented by sentences in the larger
language saying that the constants represent the two endpoints, it becomes model complete.

Similar examples, requiring arbitrarily many constants, are provided by the theory of the
usual linear order on $\A_n = \cup_{i=0}^n [2n,2n+1]_{\Q}$.  Here $(2n+2)$ constants
are necessary, to name the end points of all the different intervals.  The proof is much the same
as for a single interval.

One naturally asks if the theory of the linear order
$\mathcal Z=\cup_{n\in\Z}[2n,2n+1]_{\Q}$ (under the usual $<$ on $\Q$) becomes model complete
when the language is augmented by infinitely many constants.  Here, of course, the new constants
$c_n$ and $d_n$ (for all $n\in\Z$) should be used to name the elements $2n+1$ and $2n+2$.
The theory of $\mathcal Z^*$, with these constants, specifically describes each open interval 
$(d_{n-1},c_n)$ as dense without endpoints, and each $(c_n,d_n)$ as empty.
However, the theory says nothing about the existence of elements $x$ satisfying $d_n<x$
for every $n$, nor about elements to the left of all successor pairs:  one could place another copy of $\mathcal Z$ (without any constants), or several copies of
$\mathcal Z$, and the lack of constants allows one to show that this theory is not model complete.
Similar problems arise even when the (infinitely many) successor pairs named by constants
are arranged in different orders, as the order on the pairs of constants will always
contain either a copy of $\omega$ or a copy of its reverse order $\omega^*$,
and the theory will be unable to specify what sits at the right end of the copy of $\omega$,
or at the left end of the copy of $\omega^*$.

A more intriguing example arises when we consider the linear order $\B$ given by
the lexicographic order on $\Q\times\{0,1\}$.  This arises often in the literature, either as
the \emph{shuffle sum} of countably many copies of the two-element order, or as ``$\Q$
with every point doubled,'' or other ways.  The theory $T$ of $\B$ is not model complete:
the substructure
$$ \A=\set{(q,k)\in\B}{q<4\text{~or~}q>5}\cup\{ (4,0),(5,1)\}$$
of $\B$ is not an elementary substructure, yet satisfies $T$, being isomorphic to $\B$.
To make this $T$ model complete, we adjoin a binary relation symbol $\Adj$ to the language,
with axioms stating that $\Adj(x,y)$ holds just if $x<y$ and $\forall z\neg (x<z<y)$.
The theory $T^*$ of the extension $\B^*$ in the new language (with $\Adj^{\B^*}$ defined
as instructed) includes sentences saying that every element is half of a unique adjacency, and
that the pairs satisfying $\Adj$ are dense with no least or greatest such pair.  One then proves
that every substructure of $\B^*$ modelling $T^*$ must be elementary, so that $T^*$
is indeed model complete.

It is natural to ask whether the same could have been accomplished by augmenting the original
signature $\{ <\}$ by countably many constants, say $c_q$ and $d_q$ for all $q\in\Q$,
using $c_q$ and $d_q$ to name the elements $(q,0)$ and $(q,1)$ in $\B$, respectively.
If $\B^c$ is the structure $\B$ thus enriched, and $T^c$ is its theory, then $\B^c$ has
no substructure except itself.  However, consider the structure $\A^c$ of the lexicographic order
on $\Q\times\{0,1\}$ in which $c_q^{\A}=(q,0)$ for $q\leq 3$, but $c_q^{\A}=(3+q,0)$ for $q>3$
(and each $d_q^{\A}$ is the immediate successor of $c_q^{\A}$, necessarily).  This is still a model of $T^c$,
but in the gap between $3$ and $6$, where no constants were used, we have the same problem
as before:  the substructure $\C^c\subseteq\A^c$ containing those elements of $\A^c$ which are
$\leq (4,0)$ and those which are $\geq (5,1)$ is another model of $T^c$, but not an elementary
substructure of $\A^c$.  So in this situation, constants do not suffice; the new relation symbol
$\Adj$ was necessary to yield model completeness.

\section{Relatively Decidable Isomorphism Types}
\label{sec:thms}

Our ultimate goal is to characterize those theories that are
relatively decidable -- that is, those theories $T$ such that every
model of $T$ with domain $\omega$ is relatively decidable.
This question remains open, and we discuss it in Section
\ref{sec:precomplete}.  As a step towards that goal, we prove
the analogous result here for isomorphism types,
characterizing those structures $\A$ such that every
copy of $\A$ is relatively decidable.  In Corollary
\ref{cor:low} we will derive a near-effectiveness result
regarding those $\A$ for which not every copy is relatively decidable.

\begin{thm}
\label{thm:URD}
Let $\A$ be a countable structure in a finite relational language,
and assume that every structure (on the domain $\omega$) isomorphic
to $\A$ is relatively decidable.  Then there is a finite tuple 
$\vec{a} \in A$ such that $(\mathcal{A}, \vec{a})$ is uniformly relatively decidable. 
\end{thm}
The restriction requiring the language to be finite and relational will be removed in Corollary
\ref{cor:anylang}.
\begin{pf}
We assume there is no tuple $\vec{a}$ such that $(\mathcal{A}, \vec{a})$ is 
uniformly relatively decidable, and construct $\mathcal{B} \cong \mathcal{A}$
such that $\Phi_e^{\Delta(\mathcal{B})} \neq E(\mathcal{B})$ for every $e$.

We build $\mathcal{B}$ by constructing a generic permutation $g$ of $\omega$ and defining
$\mathcal{B}$ with domain $\omega$ such that $g: \mathcal{B} \rightarrow 
\mathcal{A}$ is an isomorphism. We follow the method in Ash, Knight, Manasse and Slaman
\cite{AKMS89}, except we work with finitary formulas rather than in a countable fragment of 
$\mathcal{L}_{\omega_1,\omega}$.  In fact, Corollary \ref{cor:low} will show
that full-scale forcing is hardly necessary to prove this theorem,
but the method is familiar to many readers and will be readily understood.

A condition is a finite partial 1-to-1 function from an initial segment of $\omega$ into $\omega$.
We define extension by $q \leq p$ if and only if 
$q \supseteq p$. We expand the language $L$ of $\A$ to $L^f$ by adding a function symbol
$f$ to denote the generic. We let $L^f(\mathcal{A})$ denote the further extension 
when we add constants for the elements of the domain $\omega$ of $\mathcal{A}$.
For a sentence $\varphi \in L^f(\mathcal{A})$ and a condition $p$, we define 
$p \Vdash \varphi$ as in \cite{AKMS89}. 

Since it is dense to add new elements to the domain and range of a condition, the resulting generic will be a permutation of $\omega$. Given a permutation $g$, 
we define the corresponding model $\mathcal{B}_g$ by 
\[
\mathcal{B}_g \models \varphi(b_1, \ldots, b_k) \Leftrightarrow \mathcal{A} \models \varphi(g(b_1), \ldots, g(b_k))
\]
for all atomic formulas $\varphi(\vec{x})$ and tuples $\vec{b}$. We use $g(\vec{b})$ to denote the tuple $\langle g(b_1), \ldots, g(b_k) \rangle$. 
Note that by definition, $g: \mathcal{B}_g \rightarrow \mathcal{A}$ is an isomorphism. 

In a similar way, each condition $p$ determines a finite part of the atomic diagram of an isomorphic copy of $\mathcal{A}$. We use $\Delta(p)$ to denote 
the finite set of atomic and negated atomic sentences already determined by $p$. For a Turing functional $\Phi$, we write $\Phi^{\Delta(p)}(n) \downarrow$ if 
the computation with finite oracle $\Delta(p)$ halts and only queries atomic and negated atomic facts determined by $p$.

For a condition $p$, we let $\vec{d}_p$ denote the domain of $p$ and $\vec{a}_p$ denote the range $p$. We say that a structure $\mathcal{C}$ 
extends $p$ if there is an isomorphism $h: \mathcal{C} \rightarrow \mathcal{A}$ such that $h \upharpoonright \vec{d}_p = p$. Note that $h$ is an isomorphism 
from $(\mathcal{C}, \vec{d}_p)$ onto $(\mathcal{A}, \vec{a}_p)$ and that if $\Phi^{\Delta(p)}(n) \downarrow$ and $\mathcal{C}$ extends $p$, then 
$\Phi^{\Delta(\mathcal{C})}(n) = \Phi^{\Delta(p)}(n)$.

Let $p$ be a condition and let $(\mathcal{E}, \vec{e}) \cong (\mathcal{A}, \vec{a}_p)$. By permuting the domain of $\mathcal{E}$, there is a structure $\mathcal{C}$ which 
extends $p$ and satisfies $(\mathcal{C}, \vec{d}_p) \cong (\mathcal{E}, \vec{e})$. That is, we can transform any isomorphic copy of $ (\mathcal{A}, \vec{a}_p)$ into 
a copy of $\mathcal{A}$ which extends $p$. Moreover, we can find $(\mathcal{C}, \vec{d}_p)$ uniformly in $(\mathcal{E}, \vec{e})$ and $p$. 

We need to build our generic $g$ so that $\mathcal{B}_g$ satisfies  
\[
\R_e: \Phi_e^{\Delta(\mathcal{B}_g)} \neq E(\mathcal{B}_g)
\]
for each $e$. We say that a condition $p$ satisfies $\R_e$ if one of the following two conditions holds.
\begin{enumerate}
\item[(C1)] There is an $n$ such that for every $\mathcal{C}$ extending $p$, $\Phi_e^{\Delta(\mathcal{C})}(n) \uparrow$. 
\item[(C2)] There is a formula $\varphi(\vec{x})$ such that either 
\[
\Phi_e^{\Delta(p)}(\ulcorner \varphi(\vec{d}_p) \urcorner) = 0 \, \text{ and } \, p \Vdash \varphi(\vec{d}_p)
\]
or 
\[
\Phi_e^{\Delta(p)}(\ulcorner \varphi(\vec{d}_p) \urcorner) = 1 \, \text{ and } \, p \Vdash \neg \varphi(\vec{d}_p)
\]
\end{enumerate} 

To see why satisfying (C1) or (C2) is sufficient, let $g$ be a generic extending $p$. If (C1) holds for $p$, then $\Phi_e^{\mathcal{B}_g}$ is not total and hence satisfies 
$\R_e$. If (C2) holds for $p$, then $\Phi_e^{\Delta(\mathcal{B}_g)}(\ulcorner \varphi(\vec{d}_p) \urcorner) = 0$ if and only if $\varphi(\vec{d}_p) \in E(\mathcal{B}_g)$, and 
hence $\mathcal{B}_g$ satisfies $\R_e$. 

To complete our theorem, it suffices to show that the set of conditions satisfying each $\R_e$ is dense.
Fixing a condition $p$, we now show that there exists an extension of 
$p$ which satisfies $\R_e$, by considering three different cases and
then showing that one of these three cases must hold. 

\medskip
\textbf{Case 1.} There is a structure $\mathcal{C}$ extending $p$, a formula $\varphi(\vec{y},\vec{z})$ and a tuple $\vec{c} \in \mathcal{C}$ disjoint from 
$\vec{d}_p$ such that 
\[
\Phi_e^{\Delta(\mathcal{C})}(\ulcorner \varphi(\vec{d}_p, \vec{c}) \urcorner) = 0 \, \text{ and } \, \mathcal{C} \models \varphi(\vec{d}_p, \vec{c}).
\]

In this case, fix an isomorphism $h: \mathcal{C} \rightarrow \mathcal{A}$ such that $h \upharpoonright \vec{d}_p = p$. Let $q$ be a condition extending $p$ such that the 
elements of $\vec{c}$ are in the domain of $q$ and $q(\vec{c}) = h(\vec{c})$, and furthermore, $q$ agrees with $h$ on a long enough initial segment of $\omega$ (i.e.~the 
domain of $\mathcal{C}$) that every atomic or negated atomic fact queried in the computation from $\Delta(\mathcal{C})$ is determined by $q$. 
Thus, we have that $q \leq p$, $\mathcal{C}$ is an extension of 
$q$, and $\Phi_e^{\Delta(q)}(\ulcorner \varphi(\vec{d}_p, \vec{c}) \urcorner) = 0$. Since the elements of $\vec{d}_p$ and $\vec{c}$ are contained in $\vec{d}_q$, 
we can view $\varphi(\vec{y},\vec{z})$ as a formula $\varphi(\vec{x})$ such that $\varphi(\vec{d}_p, \vec{c})$ is identical to $\varphi(\vec{d}_q)$. In this notation, we 
have $\mathcal{C} \models \varphi(\vec{d}_q)$ and $\Phi_e^{\Delta(q)}(\ulcorner \varphi(\vec{d}_q) \urcorner) = 0$. 

By the definition of $q$, the isomorphism $h: \mathcal{C} \rightarrow \mathcal{A}$ satisfies $h \upharpoonright \vec{d}_q = q$. Therefore, since $\mathcal{C} \models 
\varphi(\vec{d}_q)$, we have $\mathcal{A} \models \varphi(\vec{a}_q)$. Let $g$ be a generic extending $q$. By the definition of $\mathcal{B}_g$, $g$ is an 
isomorphism from $\mathcal{B}_g$ to $\mathcal{A}$ such that $g \upharpoonright \vec{d}_q = q$. Therefore, since $\mathcal{A} \models \varphi(\vec{a}_q)$,  
we have $\mathcal{B}_g \models \varphi(\vec{d}_q)$.

$\mathcal{B}_g$ is a generic structure, so the fact that $\mathcal{B}_g \models \varphi(\vec{d}_q)$ must be forced. Fix a condition $r \leq q$ such that 
$r \Vdash \varphi(\vec{d}_q)$. Because $r \leq q$, we maintain $\Phi_e^{\Delta(r)}(\ulcorner \varphi(\vec{d}_q) \urcorner) = 0$, and therefore $r$ is the 
desired extension of $p$ satisfying $\R_e$ through condition (C2). 

\medskip
\textbf{Case 2.} There is a structure $\mathcal{C}$ extending $p$, a formula $\varphi(\vec{y},\vec{z})$ and a tuple $\vec{c} \in \mathcal{C}$ (disjoint from 
$\vec{d}_p$) such that 
\[
\Phi_e^{\Delta(\mathcal{C})}(\ulcorner \varphi(\vec{d}_p, \vec{c}) \urcorner) = 1 \, \text{ and } \, \mathcal{C} \models \neg \varphi(\vec{d}_p, \vec{c}).
\]
This case proceeds as in the previous case except we work with the formula $\neg \varphi$ in place of $\varphi$. 

\medskip
\textbf{Case 3.} There is a condition $q \leq p$ and an $n$ such that 
$\Phi_e^{\Delta(\mathcal{C})}(n) \uparrow$ for all structures $\mathcal{C}$ extending $q$. 
In this case, $q$ is already an extension of $p$ satisfying $\R_e$ through condition (C1). 

\medskip

To finish the proof, we need to show that we must be in one of the three cases above.
Assume that none of the three cases apply for some fixed condition $p$ with $\R_e$. 
\comment{
This is the point were I run into trouble. I think we want to show that the structure $(\mathcal{A}, \vec{a}_p)$ is uniformly relatively decidable for a contradiction. I will first lay out 
what I think we can conclude and then say exactly where I'm stuck in trying to finish the proof.
}
Under this assumption, we now describe a uniform procedure that, given a structure
$(\mathcal{E}, \vec{e}) \cong (\mathcal{A}, \vec{a}_p)$, a formula $\varphi(\vec{x}, \vec{e})$ and a 
tuple $\vec{u}$ from $\mathcal{E}$, determines whether $\mathcal{E} \models \varphi(\vec{u}, \vec{e})$.
(Without loss of generality, $\vec{e}$ and $\vec{u}$ are disjoint.)
So, fix $(\mathcal{E}, \vec{e})$, $\varphi(\vec{x}, \vec{e})$ and $\vec{u}$. As noted above, we can use permutations of finite initial segments of $\omega$ to 
transform $(\mathcal{E}, \vec{e})$ into $(\mathcal{C}, \vec{d}_p)$ such that $(\mathcal{E}, \vec{e}) \cong (\mathcal{C}, \vec{d}_p)$ and $\mathcal{C}$ extends $p$. 
Let $\vec{c}$ be the image of $\vec{u}$ under the appropriate permutation. We have reduced our question to the following: describe a procedure that, given a structure 
$\mathcal{C}$ extending $p$, a formula $\varphi(\vec{x}, \vec{d}_p)$ and a tuple $\vec{c}$ from $\mathcal{C}$, determines, uniformly in an oracle for $\Delta(\C)$, whether $\mathcal{C} \models \varphi(\vec{c}, \vec{d}_p)$. 

For simplicity, we assume that, for each $n$, the G\"odel coding of atomic sentences
about elements of the domain $\omega$ numbers all sentences about $\{ 0,\ldots,n\}$
before it comes to any sentence about $(n+1)$.  
We then define the numbers $l_n$ (effectively) so that
the restriction $\Delta(\C)\res\{0,\ldots,n\}$ of the atomic diagram of a structure $\C$
to a finite initial segment $\{0,\ldots,n\}$ of its domain is a string in $2^{l_n}$

As described above, we may start with a structure $(\C,\dvec_p)\cong (\A,\avec_p)$.
Given a formula $\phi(\xvec)$ and a tuple
$\cvec$ from $\C$, we determine whether $\C\models\phi(\cvec)$ by searching for the following,
using our $\Delta(\C)$-oracle:
\begin{itemize}
\item
an $n\geq\max(\cvec,\dvec_p)$ and a $\sigma\in 2^{l_n}$ such that $\sigma$ extends the
characteristic function of $\Delta(\C)\res (\cvec,\dvec_p)$
and $\Phi_e^\sigma(\ulcorner \varphi(\vec{c}, \vec{d}_p) \urcorner)\converges$; and
\item
a tuple $\bvec\in\C^n$ of distinct domain elements of $\C$ which includes all of $\cvec$ and all of $\dvec_p$,
such that the map $\rho(i)=b_i$ fixes $\cvec$ and $\dvec_p$ pointwise, and such that applying
the map $\rho$ to the elements $\{0,\ldots,n\}$ sends $\sigma$ to $\Delta(\C)\res\bvec$.
\end{itemize}
We argue below that such an $n$, $\sigma$ and $\bvec$ must exist, and that
for every such collection of elements, we must have $\Phi_e^{\sigma}(\ulcorner \varphi(\vec{c}, \vec{d}_p) \urcorner)$
equal to $1$ if $\C\models\varphi(\cvec,\dvec_p)$ and equal to $0$ if not.  First, we give
the intuition for this search.

Each $\sigma$ that we find represents some atomic diagram on a set of $n$ elements, in the given language.
(Of course, it is important for this language to be finite and relational; otherwise such a $\sigma$ might
need to be infinitely long.)  Now $\Phi_e$ may converge using oracles $\sigma$ which have nothing
to do with $\C$ or with $\A$.  However, if we also find a tuple $\bvec$ in $\C$ and a bijection $\rho$
as described, then we know that the configuration described by $\sigma$ does actually occur
(with the elements $\bvec$, not necessarily with $\{0,\ldots,n\}$) in $\C$.
Moreover, it occurs with the particular $\cvec$ and $\dvec_p$ that
matter to us, since $\rho$ fixes these elements.  If we extend $\rho$ to a permutation $f$
of $\omega$, then we get a structure $(\D,\cvec,\dvec_p)$ on the domain
$\omega$, defined so that $f$ is an isomorphism from this structure onto $(\C,\cvec,\dvec_p)$.
Thanks to our specifications regarding $\rho$, the atomic diagram $\Delta(\D)$ must restrict to $\sigma$.
Therefore, $\Phi_e^{\Delta(\D)}(\ulcorner \varphi(\vec{c}, \vec{d}_p) \urcorner)$ converges to the
value $\Phi_e^\sigma(\ulcorner \varphi(\vec{c}, \vec{d}_p) \urcorner)$ that we found.
But by the failure of Cases 1 and 2, this value must be ``correct,'' in $\D$, i.e.,
must tell us accurately whether $\D\models\varphi(\cvec,\dvec_p)$.
Since the isomorphism $f$ from $\D$ onto $\C$ fixes $\cvec$ and $\dvec_p$, it also
tells us whether $\C\models\varphi(\cvec,\dvec_p)$, which is what we wanted to know.
The last point is that some such $n$, $\sigma$ and $\bvec$ do exist:  if not, then Case 3
would not have failed.  Therefore, the search described above must eventually find such
an $n$, $\sigma$ and $\bvec$, and from them we can decide whether $\C\models\varphi(\cvec,\dvec_p)$
or not, proving the uniform relative decidability of $(\A,\dvec_p)$, contrary to our original hypothesis.

Now for the full details.
We know $(\C,\dvec_p)\cong(\A,\avec_p)$, so fix an isomorphism $f:(\C,\dvec_p)\to\A(\avec_p)$.
Let $q\subseteq f$ be a partial function on an initial segment of $\omega$
long enough that $(\cvec,\dvec_p)\subseteq\dom{q}$.
Let $\avec=q(\cvec)$ in $\A$ and note that $q$ extends $p$.

Since Case 3 fails, there is some structure $\D$,
isomorphic to $\A$ via an extension of $q$, such that
$\Phi_e^{\Delta(\D)}(\ulcorner \varphi(\vec{c}, \vec{d}_p) \urcorner)\converges$.
Indeed, we may take this $\D$ to be generic, since only a finite
initial segment of $\Delta(\D)$ is actually used in this computation.
By genericity, and since Cases 1 and 2 fail,
the output of the computation is the correct answer about
whether $\D\models\varphi(\cvec,\dvec_p)$.  Fix a sufficiently large $n$
that, with $\sigma=\Delta(\D)\res l_n$, we have
$\Phi_e^{\sigma}(\ulcorner \varphi(\vec{c}, \vec{d}_p) \urcorner)\converges$.  Now $\C$
and $\D$ are both isomorphic to $\A$ via isomorphisms extending $q$, so we can take
$\rho:\{0,\ldots,n\}\to\C$ to be the restriction of the resulting isomorphism from $\D$ onto $\C$.
The elements $\bvec$ in the image of this $\rho$ are the necessary tuple from $\C$, since $\rho$,
being the restriction of an isomorphism, must send $\sigma$ to $\Delta(\C)\res\bvec$.
This proves that there do exist an $n$, $\sigma$, and $\bvec$ as described in the program
for deciding the elementary diagram of $\C$, so our search there must eventually terminate.

Of course, the $n$, $\sigma$, and $\bvec$ found by the program's search are not necessarily
the ones determined above using the isomorphism $f$.  We now prove that, for \emph{every}
$n$, $\sigma$, and $\bvec$ satisfying the conditions in the search, the value of
$\Phi_e^{\sigma}(\ulcorner \varphi(\vec{c}, \vec{d}_p) \urcorner)$ correctly describes whether
$\C\models\varphi(\cvec,\dvec_p)$ or not.  The argument was summarized when we
gave the intuition for the proof.  Define the map $\rho:\{0,\ldots,n\}\to\bvec_i$ sending $i$
to the coordinate $b_i$ in $\bvec$.  By the conditions in the search, $\rho$ is the identity map
on $\cvec$ and on $\dvec_p$, and replacing each $i$ by $\rho(i)$ in the formulas
named by the G\"odel numbers $\leq l_n$ converts $\sigma$ into $\Delta(\C)\res\bvec$.
Therefore, $\rho$ extends to an isomorphism $h$ from some structure $\D$ onto $\C$,
where $\sigma$ is an initial segment of $\Delta(\D)$, and we may assume this
$\D$ to be generic, since only a finite initial segment of $\Delta(\D)$
is prescribed.  Moreover, $(\D,\cvec,\dvec_p)\cong
(\C,\cvec,\dvec_p)$ via this $h$, by our conditions on $\rho$.  It follows that
$\Phi_e^{\Delta(\D)}(\ulcorner \varphi(\vec{c}, \vec{d}_p) \urcorner)\converges$,
with at most the initial segment $\sigma$ of the oracle being used in the computation.
If it outputs $0$ and $\D\models\varphi(\cvec,\dvec_p)$, then the generic structure $\D$,
the formula $\varphi(\yvec,\zvec)$, and the tuple $\cvec$ would have shown that Case 1
holds, contrary to our assumption.  Likewise, if it outputs $1$ and
$\D\models\neg\varphi(\cvec,\dvec_p)$, then Case 2 would have held.
Therefore, the output of $\Phi_e^{\sigma}(\ulcorner \varphi(\vec{c}, \vec{d}_p) \urcorner)$
correctly describes whether $\D\models\varphi(\cvec,\dvec_p)$.  But the isomorphism
$h$ shows that $\D\models\varphi(\cvec,\dvec_p)$ if and only if $\C\models\varphi(\cvec,\dvec_p)$,
so the conclusion of our program does correctly decide the truth of $\varphi(\cvec,\dvec_p)$
in $\C$, using only the atomic diagram of $(\C,\dvec_p)$ as an oracle.
Therefore, under the assumption that all three cases fail at $p$, $\A$ would indeed
have a uniformly relatively decidable expansion by constants.
\qed\end{pf}

We now remove the restriction to finite relational languages.

\begin{cor}
\label{cor:anylang}
Let $\A$ be a countable structure in a computable language,
such that all copies of $\A$ on the domain $\omega$ are relatively decidable.
Then there is a finite tuple 
$\vec{a} \in A$ such that $(\mathcal{A}, \vec{a})$ is uniformly relatively decidable. 
\end{cor}
\begin{pf}
To apply Theorem \ref{thm:URD}, we need the language to be finite and relational.
This is a straightforward application of the well-known theorem of Hirschfeldt,
Khoussainov, Shore, and Slinko in Appendix A of \cite{HKSS02}, which takes an arbitrary
automorphically nontrivial structure in an arbitrary computable language
and produces a directed graph with exactly the same
computable-model-theoretic properties.  Applying this to the language
of $\A$ gives a symmetric irreflexive graph $G$ in the language with
equality and a single binary relation symbol, to which we can apply
Theorem \ref{thm:URD}.

To explain more fully:  \cite{HKSS02} provides an \emph{effective
bi-interpretation} between $\A$ and $G$, in the sense of \cite[\S 5.1]{M14}.
The construction can also be viewed as a use of \emph{computable
functors} between the category of presentations of $\A$ and the category
of presentations of $G$, as defined in \cite{MPSS18}; this is shown in \cite{HTM3}
to be equivalent to effective bi-interpretation.  Using the effective interpretations,
we can build from each copy $\A'$ of $\A$ a copy $G'$ of $G$, with
$\Delta(\A')\equiv_T\Delta(G')$, and conversely.  If follows that every copy $G'$
of $G$ is relatively decidable:  given any sentence $\alpha(\gvec)$
about elements $\gvec$ of $G'$, we can translate this into a sentence
about finitely many elements $\avec$ in the copy $\A'$ of $\A$ built from $G'$,
and then determine the truth of the sentence in $\A'$
by applying the procedure for deciding $E(\A')$ from $\Delta(\A')$
(using $\Delta(G')$ to decide $\Delta(\A')$, since the interpretation is
effective).  Thus every $G'\cong G$ is relatively decidable,
and Theorem \ref{thm:URD} gives a finite tuple $\gvec$ in $G$ such that
$(G,\gvec)$ is uniformly relatively decidable.  This tuple $\gvec$
corresponds to a finite tuple $\avec$ from $\A$ (namely, the elements
of the tuples in $\A$ that interpret the elements of $\gvec$), and the uniform
decision procedure for copies $(\A',\avec~\!')$ of $(\A,\avec)$ begins by building
the structure $G'$ corresponding to $\A'$ and identifying the tuple $\gvec~\!\!'$
in $G'$ (using $\avec~\!'$), which can be done uniformly using the effective
interpretation.  Then translate the given sentence about $(\A',\avec~\!')$
into an equivalent one about $(G',\gvec~\!')$ and apply the uniform decision procedure
for copies of $(G,\gvec)$ to determine its truth in $(G',\gvec~\!')$.

It remains to consider the automorphically trivial models $\A$ of $T$.
Unsurprisingly, this is simple.
By definition, the condition means that there is a finite tuple $\avec$
from $\A$ such that every permutation of $\omega$ which fixes $\avec$ pointwise
is an automorphism of $\A$.  In this case we take that $\avec$ as our tuple
of constants.  Now for every $(\C,\cvec)\cong(\A,\avec)$,
we can find an isomorphism from $\A$ onto $\C$, effectively in
$\Delta(\C,\cvec)$:  just map each $a_i$ to $c_i$, extend this
to a permutation of a finite initial segment of $\omega$, and then extend
by the identity map to a permutation $f_{\C}$ of all of $\omega$.
This $f_{\C}$ is an isomorphism from $\A$ onto $\C$,
computable uniformly in the atomic diagram (in fact,
$\cvec$ is all that is needed), and so we may use it to decide
the truth of formulas $\varphi(\cvec,\dvec)$ in $\C$ by determining
the truth of $\varphi(\avec,f(\dvec))$ in the relatively decidable
structure $\A$.  (The oracle $\Delta(\A)$ for running the decision
procedure for $\A$ is at hand, since we know both $\Delta(\C)$
and the isomorphism $f$.)
\qed\end{pf}

The proof in Theorem \ref{thm:URD} lends itself to an effective
construction establishing the following corollary.
\begin{cor}
\label{cor:low}
Suppose that $\A$ is a countable structure on the domain $\omega$
such that, for every tuple $\avec$ of constants from $\A$, $(\A,\avec)$
fails to be uniformly relatively decidable.  Then there exists a structure
$\B\cong\A$ which is low relative to $\A$ and is not relatively decidable:
it satisfies $(\Delta(\B))'\leq_T (\Delta(\A))'$ and $E(\B)\not\leq_T\Delta(\B)$. 
\end{cor}
\begin{pf}
This is a finite-injury construction below a $\Delta(\A)$-oracle,
using the result (established
in the proof of Theorem \ref{thm:URD}) that, for an $\A$ which is not
uniformly relatively decidable, the three Cases cannot all fail.
Of course, we may assume that $\A$ itself is relatively decidable,
since otherwise the corollary is trivial.  The goal is to construct
a permutation $f$ of $\omega$ such that we can pull $\A$ back
via $f$ to form the desired structure $\B$ with $f:\B\to\A$ an isomorphism.
At each stage $s$ of the construction, we will compute (using $\Delta(\A)$) a finite
injective portion $p_s$ of $f$.  These maps $p_s$ will not all be compatible,
but with $f=\lim_s p_s$ we will satisfy the theorem.  (If the $p_s$ were
all compatible, then $f$ would be computable from $\Delta(\A)$,
and so would $\Delta(\B)$.  We conjecture that there exist
structures $\A$ for which all copies $\B$ with $\Delta(\B)\leq_T\Delta(\A)$
are relatively decidable, in which case this would be impossible.)

We mix the requirements $\R_e$ from Theorem \ref{thm:URD}
with standard lowness requirements for $f$ relative to $\Delta(\A)$:
$$ \L_e:~~(\exists^\infty s~\Phi_{e,s}^{p_s}(e)\converges)
\implies \Phi_e^f(e)\converges,$$
and with surjectivity requirements $\mathcal S_y$, each ensuring that
the element $y$ in the domain $\omega$ of $\A$ lies
in the image of all maps $p_s$ with $s\geq y$.
We fix a computable ordering $\prec$, in order type $\omega$,
of the injective elements of $\omega^{<\omega}$.

Naturally $p_0$ is the empty map.  At stage $s+1$, we have $p_s$
from the preceding stage, and we build $p_{s+1}$ by going through
the requirements in order, first $\L_0$, then $\mathcal S_0$, then $\R_0$,
then $\L_1$, and so on up to $\R_s$.  We start with $p$ as the empty string
and extend it as we encounter each of these requirements;
$p_{s+1}$ will be the ultimate result.  For a requirement $\L_e$,
we check whether $\Phi_{e,s}^q(e)$ converges for any of the
first (under $\prec$) $s$-many injective extensions
$q\leq p$ in $\omega^{<\omega}$.  If so, then we extend $p$ to the $\prec$-least
such $q$; if not, then $p$ stays unchanged.  For a requirement $\mathcal S_y$,
we leave $p$ unchanged if $y\in\rg{p}$, or else extend $p$
to map the next number $|p|$ to $y$, thus satisfying $\mathcal S_y$.

For a requirement $\R_e$, we search through the first (under $\prec$) $s$ extensions $q\leq p$
in $\omega^{<\omega}$, and through all formulas $\alpha(\xvec)$ and all $\bvec\in(\dom{q})^{<\omega}$
with $\ulcorner\alpha(\bvec)\urcorner\leq s$, checking each to see
whether $\Phi_e^{\Delta(q^{-1}(\A))}(\ulcorner\alpha(\bvec)\urcorner)\converges$.
(Here $q^{-1}(\A)$ is the finite relational structure we get by pulling back
the structure $\A$ via the finite map $q$.)  Since we have an oracle for $\Delta(\A)$,
this is computable.  The oracle also allows us to compute $E(\A)$, since by assumption
$\A$ itself is relatively decidable, and so, for each $q$ where convergence occurs,
we may check whether it agrees with the truth in $\A$ of $\alpha(q(\bvec))$.
If we find any $q$ here for which either
\begin{align*}
&[\Phi_e^{\Delta(q^{-1}(\A))}(\ulcorner\alpha(\bvec)\urcorner)\converges=0
\text{~~and~~}\A\models\alpha(q(\bvec))],\\
\text{or~~}&
 [\Phi_e^{\Delta(q^{-1}(\A))}(\ulcorner\alpha(\bvec)\urcorner)\converges=1
\text{~~and~~}\A\models\neg\alpha(q(\bvec))],
\end{align*}
then we set $p$ to equal the $q$ in the least such pair $\la q,\ulcorner\alpha(\bvec)\urcorner\ra$.
If none of the $q$ checked here have either of these properties,
then we extend $p$ to the $\prec$-least element $p'$ such that
none of the finitely many $q$ we found with
$\Phi_e^{\Delta(q^{-1}(\A))}(\ulcorner\alpha(\bvec)\urcorner)\converges$
extends $p'$.  This completes the construction.

We claim that, for every requirement $\P$, the string $p_{s,\P}$ chosen by that requirement
at stage $s$ will stabilize on some $p_{\P}$ as $s\to\infty$.
Since $p_{s,\P}$ is an initial segment of $p_{s+1}$, this will prove
that $f=\lim_s p_s$ exists.  (The $\mathcal S$-requirements ensure that
$|p_s|\geq s$ for all $s$, so $\dom{f}=\omega$.)  By induction
we assume stabilization of the string $p_{\P}$ produced by the next-higher-priority requirement.
For $\mathcal S$-requirements, the result is clear, and likewise for a requirement $\L_e$,
which will eventually find the $\prec$-least $q$ extending $p_{\P}$
with $\Phi_e^q(e)\converges$ (if any exists),
and will choose this $q$ as $p_{s+1,\L_e}$ forever after;
or else will choose $p_{\P}$ forever if there is no such $q$.
For a requirement $\R_e$, if there is any $q\leq p_{\P}$
and any $\alpha(\bvec)$
that yield a disagreement between $\Phi_e^{\Delta(q^{-1}(\A))}(\ulcorner\alpha(\bvec)\urcorner)$
and $E(\A)$, then we will eventually find the least such pair
and will choose that $q$ as $p_{s,\R_e}$ forever after.
If there is no such $q$ and $\alpha(\bvec)$, then it follows
from the proof of Theorem \ref{thm:URD} that the Case 3 described there must hold:
some $n\in\omega$ and $q\leq p$ have $\Phi_e^{\Delta(f^{-1}(\A))}(n)\diverges$
for every permutation $f$ of $\omega$ extending $q$.  In this case,
our construction will eventually light on the $\prec$-least such $q$
and will choose it as $p_{s,\R_e}$ forever after.  This completes the induction,
allowing us to fix the permutation $f=\lim_s p_s$ and to define $\B$ to be
the unique structure with domain $\omega$ from which $f$ is an isomorphism
onto $\A$.  The $\mathcal S$-requirements ensure that $f$ has image $\omega$,
and the $\L$-requirements show $f$ to be low relative to $\Delta(\A)$ in the usual way,
so that $(\Delta(\B))'\leq_T (\Delta(\A))'$.
Finally, if there were a (least) $e$ with $\Phi_e^{\Delta(\B)}=E(\B)$,
then the corresponding $\R_e$ would never have settled on a limit
$p_{\R_e}$:  once the next-higher string $p_{\P}$ had stabilized, $\R_e$ would never
have found a string giving a disagreement, but also must never have found
a $q\leq p_{\P}$ above which, for some $\alpha(\bvec)$, no convergence ever occurred.
(If there were such a string $q$, $\R_e$ would eventually have chosen it, so that $f\supseteq q$;
but then $\Phi_e^{\Delta(\B)}=E(\B)$ contradicts the claim that no convergence
above this $q$ ever occurred.)  So indeed our $\B$ fails to be relatively decidable.
\qed\end{pf}

\section{Model Precompleteness}
\label{sec:precomplete}

With the preceding results proven, it remains to ask
whether some theorem similar to Theorem \ref{thm:URD} holds
in the situation where every model of a c.e.\ theory $T$ is relatively decidable.
Each individual model of $T$ will have the property described there, of
having a uniformly relatively decidable expansion by finitely many constants,
but it is not clear whether this can be done uniformly across all the countable
models of $T$.  Here, building on the examples from Section \ref{sec:examples},
we will offer a conjecture about this situation, involving a weakened
version of model completeness (Definition \ref{defn:precomplete}).
We encourage both model theorists and computable model theorists
to examine this weaker notion, both in hope of a proof of Conjecture
\ref{conj:precomplete} and more generally to see what can be made of
it within pure model theory.

\begin{defn}
\label{defn:conservative}
For a theory $T$ in a language $\L$, let $\L^*$ be an expansion of $\L$
by finitely many constants.  An expansion $T^*$ of $T$ by new formulas
in the language $\L^*$ is a \emph{conservative expansion} if, for every
$\L$-formula $\phi$ such that $T^*\proves\phi$, we already had $T\proves\phi$.
\end{defn}

For example, the theory $T$ of $(\omega,S)$ (in the language
with just $S$ and $=$) has an expansion $T^*$ by a constant $c_0$
and a formula $\forall x~c_0\neq Sx$.  This expansion is conservative,
mainly because the obvious formula provable from $T^*$,
namely $\exists y\forall x~y\neq Sx$, was already in $T$.
Likewise, we have a conservative expansion of the theory of dense
linear orders with end points by two constants, along with formulas
stating that the two constants are the end points.

\begin{defn}
\label{defn:precomplete}
A theory $T$ in a language $\L$ is \emph{model precomplete} if there exist
finitely many constants $c_1,\ldots,c_m$ not in $\L$, and a set $T^*$ of formulas
in the expanded language $\L^*$ with these constants,
such that the theory
$$ \text{Cn}(T\cup T^*)$$
is model complete and is a conservative expansion of $T$.

We allow the tuple of constants to be empty.  That is, every model complete theory
is model precomplete.
\end{defn}

Thus both $\text{Th}(\omega,S)$ and the theory \textbf{DLO$^{++}$} of
dense linear orders with end points are model precomplete.
The same holds of the theories of disjoint unions of finitely many dense
linear orders with end points, and these theories show that it is important
to allow arbitrarily many constants.
In general, one thinks of the formulas in $T^*$ as definitions of the new constants.
The condition that $T^*$ should not yield proofs of any $\L$-formulas not already
in $T$ essentially restrains $T^*$ from adding information to the original $T$:
$T^*$ simply describes the new constants.  One might hope for $T^*$ to be a finite
set (equivalently, a single formula), but here we allow it to be infinite.

We include one further example, the theory \textbf{DLO$^{\pm}$} of dense linear
orders with exactly one end point.  This theory is incomplete, as it does not specify
whether the order has a left end point or a right end point, but only that it must
have exactly one of the two.  However, despite its incompleteness, this theory is
model precomplete:  adding a constant $c$, and a sentence stating that
$c$ is an end point, yields a model complete (though still incomplete) theory,
as the reader can check.  Alternatively, adding a constant $d$ and a sentence
making $d$ a left end point yields a complete, model complete theory, although
this theory is a nonconservative extension of \textbf{DLO$^{\pm}$}.

All the examples just mentioned are relatively decidable theories.
Indeed, it is clear that every
model precomplete theory $T$ is relatively decidable, although not necessarily
uniformly so:  adding the new constants makes the new theory model complete,
and so, in every model $\A$ of $T$ with domain $\omega$, we have an effective
procedure for deciding $E(\A)$ from $\Delta(\A)$, which is uniform apart from needing
to know the values of the constants in $\A$.  One would hope to prove the converse:
that every relatively decidable theory is model precomplete, and therefore becomes
\urd\ upon adjunction of the right constants and the appropriate properties
of those constants.

Notice that the properties of the new constants, for each of Th$(\omega,S)$,
\textbf{DLO$^{\pm}$}, and \textbf{DLO$^{++}$}, are not just existential formulas:
in these examples, each constant is defined by a universal formula.  With existential
definitions, we would not have needed to be given the constants, since each could
be identified just by searching for an element satisfying the correct existential definition.
On the other hand, Definition \ref{defn:conservative} does not require that the
new formulas in $T^*$ define the constants uniquely; $T^*$
simply states certain properties of the constants.  In our examples, however,
the constants are the unique elements realizing these properties.

\begin{conj}
\label{conj:precomplete}
A theory is model precomplete if and only if it is relatively decidable.
\end{conj}
As noted above, the forward direction is immediate, and all relatively
decidable theories known to us are model precomplete.  Nevertheless,
proofs of this conjecture have been elusive, and we must leave it as
an open question here.  As noted in Section \ref{sec:thms}, Theorem
\ref{thm:URD} is a good first step, showing that at the level of isomorphism
types, a phenomenon very similar to model completeness is indeed
equivalent to relative decidability.  It would be natural to try to extend
this to a proof of the conjecture for complete theories (thus replacing
isomorphism, from Theorem \ref{thm:URD}, by elementary equivalence),
and then to continue to theories more generally.

\parbox{4.7in}{
{\sc
\noindent
Department of Mathematics \hfill \\
\hspace*{.1in}  University of San Francisco \hfill \\
\hspace*{.2in}  2130 Fulton Street \hfill \\
\hspace*{.3in}  San Francisco, CA 94117  U.S.A. \hfill}\\
\hspace*{.045in} {\it E-mail: }
\texttt{jcchubb\at {usfca.edu} }\hfill \\
\medskip
\hspace*{.045in} {\it Webpage: }
\texttt{www.cs.usfca.edu/$\widetilde{\phantom{w}}$jcchubb}\hfill \\
}

\parbox{4.7in}{
{\sc
\noindent
Department of Mathematics \hfill \\
\hspace*{.1in}  Queens College -- C.U.N.Y. \hfill \\
\hspace*{.2in}  65-30 Kissena Blvd. \hfill \\
\hspace*{.3in}  Flushing, New York  11367 U.S.A. \hfill \\
Ph.D. Programs in Mathematics \& Computer Science \hfill \\
\hspace*{.1in}  C.U.N.Y.\ Graduate Center\hfill \\
\hspace*{.2in}  365 Fifth Avenue \hfill \\
\hspace*{.3in}  New York, New York  10016 U.S.A. \hfill}\\
\hspace*{.045in} {\it E-mail: }
\texttt{Russell.Miller\at {qc.cuny.edu} }\hfill\\
\medskip
\hspace*{.045in} {\it Webpage: }
\texttt{qcpages.qc.cuny.edu/$\widetilde{\phantom{w}}$rmiller}\hfill \\
}

\parbox{4.7in}{
{\sc
\noindent
Deaprtment of Mathematics \hfill \\
\hspace*{.1in}  University of Connecticut \hfill \\
\hspace*{.2in}  341 Mansfield Road U1009\hfill \\
\hspace*{.3in}  Storrs, CT 06269-1009  U.S.A.\hfill}\\
\hspace*{.045in} {\it E-mail: }
\texttt{solomon\at {math.uconn.edu} }\hfill \\
\medskip
\hspace*{.045in} {\it Webpage: }
\texttt{www.math.uconn.edu/$\widetilde{\phantom{w}}$solomon}\hfill \\
}

\end{document}